\documentclass[12pt]{article}
\usepackage{graphicx} % Required for inserting images
\usepackage{amsmath,amssymb,enumerate,xcolor}

\title{When a $(1,1)$-tensor generates separation of variables of  a certain metric.}
\author{Andrey Yu. Konyaev\footnote{Faculty of Mechanics and Mathematics, Moscow State University,  and  Moscow Center for Fundamental and Applied Mathematics, 119992, Moscow Russia
 \ \ \quad {\tt  maodzund@yandex.ru}}  \quad \& \quad  Jonathan M.\ Kress\footnote{School of Mathematics and Statistics, UNSW Sydney, NSW 2052, Australia}  \quad \& \quad  Vladimir S.\ Matveev\footnote{
Institut f\"ur Mathematik, Friedrich Schiller Universit\"at Jena,
07737 Jena,  Germany  \ \ \quad {\tt  vladimir.matveev@uni-jena.de}}}

\date{July  2023}

\newtheorem{theorem}{Theorem}[section]
\newtheorem{remark}{Remark}[section]
\newtheorem{definition}{Definition}[section]

\newtheorem{lemma}[theorem]{Lemma}
\newtheorem{example}{Example}[section]

\newcommand{\weg}[1]{}

\usepackage{biblatex} %Imports biblatex package
\begin{document}

\maketitle
\begin{abstract} 
By a $(1,1)$-tensor field  $K= K^i_j$, we construct an explicit   system of differential invariants  that   vanish if and only if   there (locally)  exists a metric for which  $K$ generates separation of variables. 

\end{abstract}

\section{Introduction}
We say that  a $(1,1)$-tensor field $K^i_j$ {\it generates 
 (orthogonal) separation of variables} for a metric $g$ of any signature, if the following conditions are fulfilled:
 it has $n:= \textrm{dim}(M)$ distinct eigenvalues, there exists a local coordinate system $(x^1,...,x^n)$ in which $K^i_j$ is given by a diagonal matrix $K= \operatorname{diag}(K_1,...,K_n)$,    and  the $(0,2)$-tensor field $K_{ij}:= g_{si} K^s_j$ is Killing. The later condition means that $K_{ij}$ is symmetric with respect to the lower indices (which, in particular, also implies that $g$ is diagonal) and the following 
{\it Killing equation} is fulfilled: 
\begin{equation} \label{eq:killingequation}
    K_{ij,k}+  K_{jk,i}+  K_{ki,j}=0.
\end{equation}
Above we denoted by comma the covariant derivative for the Levi-Civita connection of $g$.
Geometrically, this condition means that the quadratic in velocities function
$$
I:TM\to \mathbb{R}\, , \ \ I(\xi)=K(\xi,\xi)  
$$ 
is constant along the orbits of the geodesic flow of $g$  (we refer to $I$ as an {\it integral  } corresponding to $K$).
 
  It is known that 
the existence of such a $(1,1)$-tensor $K^i_j$   
 is equivalent to  the  (local)   existence of $n$ Killing tensors $\overset{1}K_{ij}= g_{ij},  \overset{2}K_{ij}= K_{ij}, \overset{3}K_{ij}, ..., \overset{n}K_{ij}$  such that they are linearly independent and diagonal in the coordinate system $x^1,...,x^n$, see e.g.\ \cite[Proposition 2.2]{Benenti56-93}.  The integrals corresponding to  these tensors, viewed as functions on $T^*M$ equipped with the standard symplectic structure,   Poisson commute. 
Moreover, in the coordinate system $(x^1,...,x^n)$, 
the tensors  $ \overset{1}K_{ij},..., \overset{n}K_{ij}$ have  
 the so-called {\it  St\"ackel form}  by \cite{eisenhart}. That is,   there exists  a  non-degenerate $n\times n$ matrix 
$S= (S_{ij})$  with $S_{ij}$ being a function  of the $i$-th variable $x^i$ only and such that 
 the following condition holds:  
 \begin{equation} 
 \label{eq:St_intro}
S I = P,  
\end{equation} 
where  $I$ is  an $n$-vector whose components are the integrals corresponding to $ \overset{1}K,..., \overset{n}K$, 
and   $P=  (p_1^2,p_2^2,..., p_n^2)^\top$ is the $n$-vector of the squares of  momenta.  
In this case, the Hamilton-Jacobi equation 
\begin{equation} \label{eq:HJ}
 \tfrac{1}{2} g^{ij}p_ip_j = c_1 \,   ,  \  p_i= \tfrac{\partial W}{\partial x^i}  
\end{equation} 
admits (locally)  a general solution of the form 
$$
W(x, c)= \sum_{i=1}^n W_i(x^i, c) \,  ,  \   c=(c_1,...,c_n) \, ,  \ \det\left[ \tfrac{\partial^2 W}{\partial x^i\partial c_j}\right]\ne 0,
$$
with  $W_i(x^i, c)=\pm \int^{x^i} \sqrt{   \sum_{s} S_{is}(\xi) c_s} \operatorname{d} \xi$.

Orthogonal separable coordinates have a long history and are widely used in applications, see e.g. \cite{BolsinovKonyaevMatveev-OrthogSep, BMF, book}.

\weg{
Orthogonal separable coordinates on a (pseudo-)Riemannian manifold $M$ with metric $g$, are given by a Killing tensor $K_{ij}$ of that metric for which the endomorphism field  $K^i_j$ has $n$ distinct eigenvalues and vanishing Haantjes torsion.

Given such an endomorphism $K^i_j$, it is interesting to know, what  conditions on it are locally equivalent to the existence of a metric $g$ for which $K_{ij}$ is a Killing tensor.

\begin{definition}[Separation system]
A non-degenerate metric $g_{ij}$ and $(1,1)-$tensor $K^i_j$ form a separation system $(g_{ij},K^i_j)$ if:
\begin{enumerate}
    \item $K_{ij} = g_{i\ell}K^{\ell}_j$ is Killing for some $g_{ij}$;
    \item $K^i_j$ has distinct eigenvalues and
    \item the Haantjes torsion of $K^i_j$ vanishes.
\end{enumerate}
\end{definition} }

In this paper we ask and answer the following natural question: 
 given $K^i_j$, what are the necessary and sufficient conditions for the existence of a metric $g_{ij}$ such that $K$ generates separation of variables for this metric?
\weg{Then, if such a metric exists, can an explicit formula be given for $g_{ij}$?}

Two necessary conditions can be read from the definition: $K$ must
have $n$ distinct eigenvalues, and there must exist a coordinate system in which $K$ is diagonal.
Recall that  the condition that $K^i_j$ has $n$ distinct eigenvalues is equivalent to the condition that the discriminant of the characteristic polynomial is different from  zero. For $(1,1)$-tensors with $n$ distinct eigenvalues, the local existence of a coordinate system such that  the tensor field is diagonal is equivalent to the vanishing of the Haantjes torsion by \cite{haantjes}. 

Note that both conditions, vanishing of the Haantjes torsion and nonvanishing of the discriminant, can be checked in any coordinate systems using the explicit formulas for the discriminant and for the Haantjes torsion.  

The goal and the main result of our paper is to give explicit formulas which can be calculated in any coordinate system and which ``tell'' whether a given tensor field $K^i_j$, with $n$ distinct eigenvalues at every point and vanishing Haantjes torsion, generates a separation of variables for some metric.

We obtained our  result  by employing two known results about integrable 
PDE-systems of hydrodynamic type.  Let us recall the necessary definitions. By a {\it system of 
hydrodynamic type,}  we understand the following quasilinear systems of $n$ PDEs on $n$ unknown functions $u^1,...,u^n$ of two coordinates $(t,x)$:   
\begin{equation}\label{sys:hidro}
\frac{\partial }{\partial  t}u^i = \sum_j A^i_j(u) \frac{\partial  }{\partial x}u^j,  
\end{equation}
where $A^i_j$
 is a matrix whose components depend on $u$ and with no explicit dependence on $t$ and $x$.  Actually, this matrix should be viewed as an $(1,1)$-tensor field, since if we change the unknown functions $u^i$  by a diffeomorphism $u_{new}= u_{new}(u_{old})$,  the equation remains a system of hydrodynamic type and the matrix $A$ of the system transforms by the tensorial rule for $(1,1)$-tensors. 

Let us recall necessary definitions from the theory of integrable systems of hydrodynamic type. 
 By \cite{tsarev}, 
 the system \eqref{sys:hidro} is called  semi-hamiltonian, if the matrix $A$ is diagonal with $n$ distinct eigenvalues, $A= \operatorname{diag}(A_1,...,A_n)$, 
 and in addition the following equations   are fulfilled for any $i\ne j \ne k \ne i$:    
\begin{equation} \label{eq:tsarev}
\frac{\partial }{\partial u^k}\left(\frac{1}{A_j - A_i} \frac{\partial A_i}{\partial u^j}\right)= \frac{\partial }{\partial u^j}\left(\frac{1}{A_k - A_i} \frac{\partial A_i}{\partial u^k}\right). 
\end{equation}
We slightly generalise this definition and say that the  system \eqref{sys:hidro} is 
{\it semi-hamiltonian}, if there \underline{exists}
a diffeomorphic change $u_{new}= u_{new}(u_{old})$ of unknown functions such that 
after this change the corresponding matrix $A_{new}$ is diagonal and the conditions  \eqref{eq:tsarev} are fulfilled.

Next, by \cite{rozhdenstvenskij}, the system \eqref{sys:hidro} is {\it weakly nonlinear}, if the derivative of every eigenvalue in the direction of its eigenvector is zero.  Since the system \eqref{sys:hidro} is fully determined by the $(1,1)$-tensor, we say that the (1,1)-tensor $A$ is {\it semi-hamiltonian} (resp., {\it weakly nonlinear}), if the corresponding system is such.  

The following theorem connects  $(1,1)$-tensors generating separation of variables and semi-hamiltonian weakly nonlinear systems of hydrodynamic type: 

\weg{
\begin{theorem}[Folklore: Eisenhart, Benenti]
    In this case there exists coordinates such that $K^i_j = \mbox{diag}(k_1,k_2,\ldots,k_n)$.  These coordinates are orthogonal separating coordinates in the standard sense, that is $g=\mbox{diag}(g_1,g_2,\ldots,g_n)$.
\end{theorem}}

\begin{theorem} \label{thm:2}
    Let $K^i_j$ be a $(1,1)$-tensor with  vanishing Haantjes torsion  and $n$ distinct eigenvalues. Then there exists (locally)   
     a metric $g$ for which $K$ generates separation if and only if  $K$ is 
    \begin{itemize}
        \item[(I)] weakly nonlinear and
        \item[(II)] semi-hamiltonian.
    \end{itemize}
\end{theorem}

This theorem, for which we give a short direct proof in Section \ref{sec:2},
was previously known
in a different form, see e.g.\  \cite{Ferapontov1,Ferapontov1992,FerapontovFordy}. 

By Theorem \ref{thm:2}, in order to understand whether a given $(1,1)$-tensor field with $n$ distinct eigenvalues and vanishing Haantjes torsion generates separation of variable for a certain (unknown)
metric $g$, it is sufficient to understand when such a tensor is weakly nonlinear and semi-hamiltonian. Theorem \ref{thm:3} below is a criterion for weak nonlinearity for gl-regular $(1,1)$-tensors at algebraically generic  points.

We say that a $(1,1)$-tensor is {\it gl-regular}, if at every point the geometric multiplicity (i.e., the dimension of the corresponding eigenspace) of every eigenvalue is one. A point is {\it algebraically generic}  for a $(1,1)$-tensor, if in a neighbourhood of the point the Segre  characteristic is constant (i.e., the eigenvalues may change but their number and the structure of the corresponding Jordan blocks remains the same). The other points are called {\it singular}. 
Clearly, a point is algebraically generic  for  a gl-regular $(1,1)$-tensor if and only if the number of the eigenvalues is constant in a neighbourhood of the point.

Standard application of the implicit function theorem shows that in  a neighborhood of an algebraically generic points, eigenvalues can be viewed as smooth functions (near  singular points  they may not be well defined, since their number may change, and even if defined they may be not smooth).

Note also 
that $(1,1)$-tensors with $n$ distinct eigenvalues are gl-regular; moreover,  all points are algebraically generic. 

%We say that a $(1,1)$-tensor is {\it weakly nonlinear}, if at any regular point the derivative of any eigenvalue in %the direction of the corresponding eigenvector  is zero.  

\begin{theorem} \label{thm:3} Let $K^i_j$ be a gl-regular tensor 
with characteristic polynomial $$\chi(t):=\det(t \operatorname{Id} - K)=  t^n- \sigma_1 t^{n-1} - \sigma_2 t^{n-2}-...-\sigma_n.$$ Then, it is weakly nonlinear  at algebraically generic  points if and only if for every point   we have:
\begin{equation}\label{eq:M3} 
  (K^*)^{n-1} d\sigma_1 +  (K^*)^{n-2} d\sigma_2+...+ (K^*)^{0} d\sigma_n=0.
\end{equation}

\end{theorem}
 In Theorem above, $K^*:T^*M\to T^*M$
    denotes the operator dual to $K:TM\to TM$,  that is, $K^*\alpha(\xi)= \alpha(K(\xi))$, and   $(K^*)^{m}:=\underbrace{K^*\circ K^*\circ \cdots \circ K^*}_{\textrm{ $m$ times}}$.
    
    In the matrix notation,  the equation \eqref{eq:M3} reads
    \begin{equation} \label{eq:M4}
    d\sigma_1 K^{n-1} + d\sigma_2 K^{n-2}+...+ d\sigma_n=0, 
    \end{equation}
    where $d\sigma_m= \left(\frac{\partial \sigma_m}{\partial x^1},\dots,\frac{\partial \sigma_m}{\partial x^n}  \right)$. 

Theorem \ref{thm:3} was first  formulated in \cite{Ferapontov1} (see equation (14) and proposition after it)
and was used many times since that. The proof was never published. The standard way(s) to formulate \ref{thm:3}, see e.g.\footnote{All these references deal with integrable systems of hydrodynamic type and may assume 
gl-regularity without mentioning it explicitly}  \cite[end of p. 392]{Ferapontov1},  \cite[\S  6]{Ferapontov2} or \cite[\S 2]{FerapontovVergallo} may make the impression   that Theorem
\ref{thm:3} holds for any, not necessary gl-regular (1,1)-tensors. The example below shows that  gl-regularity is necessary.

\begin{example}
    On $\mathbb{R}^2(x,y)$, 
    consider the (1,1)-tensor field  given by the matrix 
$$\begin{pmatrix} \lambda(x,y)& 0 \\ 0 & \lambda(x,y)\end{pmatrix}.$$
Though 
\eqref{eq:M3}  holds, the tensor field  is not always weakly nonlinear.    
\end{example}

We give a proof of Theorem \ref{thm:3} in Section \ref{sec:3}.

Finally,  Theorem \ref{thm:4} below answers whether a given $(1,1)$-tensor with vanishing Haantjes torsion and $n$ distinct eigenvalues is semi-hamiltonian. 

Let us introduce necessary notations. Consider two tensor fields  skew-symmetric in  lower indices: $Q$ of type $(1,q)$ and $R$ of type $(1,r)$. In \cite[\S 6]{Nijenhuis1955}, A. Nijenhuis introduced   the following tensorial operation,   sending the  pair  $Q,R$ to a tensor field of type $(1,q+r)$ skew-symmetric in  lower indices:

\begin{equation}\label{eq:10}
\begin{array}{ll} [Q^{i_1}_{j_1\dots j_q},R^{i_2}_{j_{q+1}\dots j_{q+r}}   ]^i_{j_1...j_{q+r}} &=\operatorname{Skew} \left(  Q^s_{j_1...j_q} \frac{\partial }{\partial x^s} R^i_{j_{q+1}\dots j_{q+r}} -R^s_{j_{q+1}...j_{q+r}} \frac{\partial }{\partial x^s} Q^i_{j_{1}\dots j_{q}} \right.
\\&
-  \left.q Q^i_{sj_2\dots j_q}\frac{\partial }{\partial x^{j_1}}  R^s_{j_{q+1}\dots j_{q+r}  } +r  R^i_{sj_{q+2}\dots j_{q+r}}\frac{\partial }{\partial x^{j_{q+1}}}  Q^s_{j_{1}\dots j_{q}  } 
\right), \end{array}
\end{equation}
where ``$\operatorname{Skew}$''
denotes the skew-symmetrisation (with division) with respect to the indexes $j_1,...,j_{q+r}$.
(The right hand side of \eqref{eq:10} has the correct number of indexes and is evidently  skew-symmetric in lower indices, so the nontrivial part of the result of Nijenhuis is that the formula \eqref{eq:10} defines a tensorial operation). 

Note also that for a (1,1)-tensor field $A$,   the $(1,2)$-tensor $[A,A]$ coincides with the Nijenhuis torsion (see e.g.  \cite[Definition 2.3]{Nijenhuis1}) which we denote by $N= N_{jk}^i$.

Next, given a (1,1)-tensor $A^i_j$ we construct, using   invariant tensorial operations and following \cite[\S 4]{Sharipov96}, the $(1,3)$- tensor field $P^i_{jk\ell}$:

First we construct (1,3)-tensor field $K$ by 
$$
K=3 [ [A,A], A^2]= 3 [ N, A^2]. 
$$
Then we construct $(1,3)$-tensor field $M$ by 
$$
\begin{array}{cl}M(X,Y,Z):= &  N(X, AN(Y, Z)) + N(AX, N(Y, Z)) - N(N(X, Z), AY)\\
& + N(N(X, Y), AZ) - N(X, N(AY, Z)) - N(X, N(Y, AZ)).\end{array}
$$
Next, we construct $(1,3)$-tensor field   $Q$ by  
$$
\begin{array}{cl}Q(X, Y, Z):= & K(AX, AY, Z) - K(A^2X, Y, Z) - K(X, AY, AZ) + K(AX, Y, AZ)\\
& + 4M(AX, Y, Z) - 2M(X, AY, Z) - 2M(X, Y, AZ).\end{array}
$$
Finally, we construct the
$(1,3)$-tensor field $P$ by
\begin{equation} \label{eq:last} P(X, Y, Z) = AQ(X, AY, Z) + AQ(X, Y, AZ) - A^2Q(X, Y, Z) - Q(X, AY, AZ).\end{equation}

\begin{theorem}[\cite{Sharipov96}] \label{thm:4}
Suppose $A^i_j$ has $n$ distinct eigenvalues and vanishing Haantjes torsion. Then, $A$ is semi-hamiltonian if and only if $P$ is identically zero.
\end{theorem}

As clearly written above, this theorem is due to \cite{Sharipov96}, see Theorem 7 there. Certain authors of the present  paper found it complicated to understand the arguments of the proof, so let us explain how one can easily 
verify it using standard computer algebra software, e.g., Maple. 

First, observe that the construction of $P$ is tensorial, so one can choose any coordinate system to verify $P=0$. Therefore, without loss of generality, we may assume that $A= \operatorname{diag}(A_1,...,A_n)$. 

Next, observe that components of $P^i_{jk\ell}$ are zero if $i\not \in \{ j,k,\ell\} $ 
and that for  $i\in \{ j, k, \ell\} $  the component $P^i_{jk\ell}$ is given by an expression involving $A_j, A_k, A_\ell$ and their derivatives with respect to $x^j,x^k,x^\ell$ only.  Moreover, the formula for the component  $P^i_{jk\ell}$ does not depend on the dimension.  Note also that the semi-hamiltonicity 
condition \eqref{eq:tsarev} involves $A_j, A_k, A_\ell$ and the derivatives with respect to $x^j$ and $x^k$ only, and does not depend on the dimension. 

This implies that it is sufficient to prove Theorem \ref{thm:4} in dimension 3 only, and this can be done with the help of e.g. Maple: One takes $$A= \textrm{diag}(A_1(x^1,x^2,x^3), A_2(x^1,x^2,x^3), A_3(x^1,x^2, x^3),$$ calculates, using Maple (note that \eqref{eq:10} is implemented in the Maple Tensor package), the tensor $P$, and compares its components with \eqref{eq:tsarev}.

\section{ Proof of Theorem \ref{thm:2}} \label{sec:2}

Theorem \ref{thm:2}, in one direction,  claims  that for $K^i_j$ which is locally diagonalisable, and who's all eigenvalues are different, the existence of a metric such that $K_{ij}$ is Killing implies that $K$ is weakly   nonlinear and semi-hamiltonian.  Let us prove this claim.

We work in the diagonal coordinates for $K^i_j$, that is, $K^i_j = \textrm{diag}(K_1,...,K_n)$. We assume the existence of the metric $g$ such that $K_{ij}$ is Killing and we want to show that the eigenvalues 
$K_i$ satisfy the  semi-hamiltonicity and weak nonlinearity conditions.

Since $K_{ij}$ is symmetric, the metric $g$ is   also diagonal in these coordinates so the corresponding Hamiltonian is
\[
H_g:= \frac{1}{2}\left(\varepsilon_1 \exp(g^1) p_1^2 +...+ \varepsilon_n \exp(g^n)p_n^2\right)
\] with $\varepsilon_i\in \{-1, 1\}.$ Then, 
 the integral corresponding to the Killing tensor field  is 
\[
F:= \varepsilon_1\exp(g^1)K_1 p_1^2 +...+ \varepsilon_n\exp(g^n)K_np_n^2.
\]
The Poisson bracket reads then: 
\begin{equation}\label{BracketHgF}
\begin{array}{lcl}\{H_g, F\}& = &   \sum_i \exp(2g^i) \frac{\partial K_i}{\partial x^i}p_i^3 \\ & +& \sum_{i\ne j}  \varepsilon_i \varepsilon_j \exp(g^j+g^i) \left( (K_i-K_j) \frac{\partial g^i }{\partial x^j} +  \frac{\partial K_i}{\partial x^j}\right)p_i^2 p_j.\end{array}
\end{equation}
Equating the $p_i^3$-coefficient of this polynomial in $p_1,...,p_n$ to zero gives us the weak nonlinearity condition 
$$\frac{\partial K_i}{\partial x^i}=0$$
as we want.  Next, equating the $p_i^2p_j$-coefficient with $i\ne j$  
to zero gives us the condition 

\begin{equation} \label{eq:M1} \frac{\partial }{\partial x^j} g^i  = \frac{1}{K_j - K_i} \frac{\partial K_i}{\partial x^j}  \end{equation}
implying  the semi-hamiltonicity condition 
\begin{equation} \label{eq:M2}
\frac{\partial }{\partial x^k}\left(\frac{1}{K_j - K_i} \frac{\partial K_i}{\partial x^j}\right)= \frac{\partial }{\partial x^j}\left(\frac{1}{K_k - K_i} \frac{\partial K_i}{\partial x^k}\right)
\end{equation}
in view of trivial equality  $$\frac{\partial }{\partial x^k} \frac{\partial }{\partial x^j} g^i =  \frac{\partial }{\partial x^j} \frac{\partial }{\partial x^k} g^i.$$

In order to prove Theorem \ref{thm:2} in the other direction, we view \eqref{eq:M1} as a system of PDEs on the functions $g^1,...,g^n$; the coefficients of this system are algebraic expressions  in 
$K_i$ and their  derivatives.  Our goal is to show that the system has a solution.  

First, observe that the system decouples in $n$  essentially the same subsystems: one on $g^1$, another on  $g^2$ etc. We will show the existence of the solution on the subsystem  on $g^1$ which is sufficient for our goals because  there is no essential
difference between $g^1$ and  other $g^j$.

The subsystem of the system \eqref{eq:M1} containing $g^1$ consists of $n-1$ equations: 
\begin{equation} \label{eq:M6}
\frac{\partial g^1}{\partial x^j} = \frac{1}{K_j - K_1} \frac{\partial K_1}{\partial x^j} \ \ \textrm{with $j>1$}.
\end{equation}
Since \eqref{eq:M6} does not contain the derivatives of the unknown function $g^1$ with respect to $x^1$, we may view $x^1$ as a parameter. 
Then, the system \eqref{eq:M6} is just the condition that the differential of the function 
$g^1$ with respect to the variables $x^2,...,x^n$ equals to 
\begin{equation}\label{eq:M7}
\sum_{j=2}^n \frac{1}{K_j - K_1} \frac{\partial K_1}{\partial x^j}  dx^j.   
\end{equation}

The semi-hamiltonicity condition 
\eqref{eq:M2}  implies that  \eqref{eq:M7} viewed as a local 1-form on $\mathbb{R}^{n-1}(x^2,...,x^n)$ is  closed. Then, for any initial value $g^1(x^1,0,...,0)$,  there exists a precisely one   solution. We therefore proved the existence of $g^1$ satisfying  \eqref{eq:M1}; we also see that the freedom is the choice of a  function  $\hat g^1(x^1) =g^1(x^1,0,\dots,0)$ of one variable.
Similarly,  we show    the existence of $g^2$ etc.  Theorem \ref{thm:2} is proved.

\begin{remark} The  calculations in the proof, e.g. those around the formula \eqref{BracketHgF} were of course done many times in the literature related to orthogonal separation of variables, see e.g.   \cite[proof of Lemma 1.2]{Benenti55-92}.]
\end{remark}

\section{Proof of Theorem \ref{thm:3}} \label{sec:3}

We take a gl-regular point $p$ and work in a coordinate system such that at the point $p$ the matrix $K^i_j$ has the block diagonal form 
\begin{equation}  \label{eq:M5} K= \operatorname{diag}(J_{m_1}(\lambda_1),...,J_{m_k} \lambda_k),\end{equation} 
where $J_m(\lambda) $ is the standard Jordan block of  
 dimension $m\times m$ with eigenvalue $\lambda$: 
$$J_m(\lambda)=\begin{pmatrix}
\lambda &         1& 0   & \cdots  & 0 \\
     0   & \ddots &  \ddots  & \ddots  &\vdots  \\
     \vdots    &\ddots  & \ddots & \ddots &    
        \\
    \vdots & &\ddots  & \lambda & 1     \\ 
   0 & \cdots &\cdots  & 0 & \lambda   
\end{pmatrix}
$$

We need  to show that \eqref{eq:M3}  is equivalent to     the condition $$ \frac{\partial \lambda_1}{\partial x^{1}}=0\ ,  \ \  \frac{\partial \lambda_2}{\partial x^{1+ m_1}}=0 \ , \ \  \dots \ , \ \frac{\partial \lambda_k}{\partial x^{1 + m_1+...+m_{k-1}}}=0. $$ 
This  follows immediately  from the next Lemma:

\begin{lemma} \label{lem:1}  Assume  $K$ is given by \eqref{eq:M5} at the point $p$. Then, at $p$, the left-hand side  of \eqref{eq:M3}  has the following components: the first $m_1-1$ are zero, the component number $m_1$ equals 
$$ m_1 \frac{\partial \lambda_1}{ \partial x^1 } 
 \prod_{s\ne 1} (\lambda_1- \lambda_s)^{m_s}, $$ 
 the next $m_2-1$ components are zero, the component number $m_1+ m_2$ equals 
$$ m_2 \frac{\partial \lambda_2}{ \partial x^{m_1+1} } 
 \prod_{s\ne 2} (\lambda_2- \lambda_s)^{m_s} $$ and so on.  
\end{lemma}
This Lemma can be 
verified for   small   $n$ using computer algebra software, e.g., Maple. Under the additional 
assumption that $K$ is diagonal which is sufficient for  Theorem \ref{thm:2} of our paper, it can be proved using standard ``Vandermonde'' identities of symmetric polynomials. Since  Theorem \ref{thm:3} can  be used and actually has been used in general case, we will give  a proof in the general case. 

{\bf Proof of Lemma \ref{lem:1}.}
Let us first observe that the assumptions and the statement of Lemma \ref{lem:1} uses only the form of $K$ at $p$ and the derivatives of the eigenvalues at $p$. Therefore, we may  replace $K$ by another (1,1)-tensor  coinciding with $K$ at $p$ and such that for every point 
the set of eigenvalues of  this new tensor coincides with that of $K$. In other words, without loss of generality we may and will  assume that in a small neighbourhood of $p$ in our  local coordinate system  $(x^1,...,x^n)$
the tensor $K$ has the form \eqref{eq:M5} with $\lambda_i$ being functions  of $(x^1,...,x^n)$.

Next, for any   
matrix-valued function $A(x)^i_j$ let us define  
$D(A)$ (which is a $n$-tuple of  functions)
 by $$D(A)_i:= \sum_s \tfrac{\partial }{\partial x^s} A^s_i. $$
The mapping  $D$ above is not covariant; choosing the flat connection $\nabla$ corresponding to the coordinates $x^1,...,x^n$ one can write it as $\nabla_s A^s_i$ and then one can understand $D$ as a mapping from (1,1)-tensors to (0,1)-tensors. Since we work in a fixed  coordinate system, we do not care about the geometric meaning of the mapping. It is useful though to rewrite the mapping $D$ in the following ``matrix'' form 
$$
D(A)= \left(\frac{\partial}{\partial x^1}, \dots , \frac{\partial}{\partial x^n}\right) A.
$$

 Next, consider $D(K^{n} -\sigma_1 K^{n-1}-...-\sigma_n \operatorname{Id}). $ It is zero by the Hamilton-Cayley Theorem. By the Leibniz rule it is equal to 
$$
-d\sigma_1 K^{n-1}-d\sigma_2 K^{n-2}-...-d\sigma_n\operatorname{Id} + D(K^n) - \sigma_1 D(K^{n-1})-...-\sigma_{n-1} D(K). 
$$
Thus, 
$$
d\sigma_1 K^{n-1}+d\sigma_2 K^{n-2}+...+d\sigma_n\operatorname{Id} = D(K^n) - \sigma_1 D(K^{n-1})-...-\sigma_{n-1} D(K). 
$$

Next, recall  that for the Jordan block $J_m(\lambda)$
 and for any polynomial $f(t)$ we have  
 $$f(J_m(\lambda)) = \begin{pmatrix}
f(\lambda) &         f'(\lambda)&   \frac{f''(\lambda)}{2}   & \cdots  & \frac{f^{(m-1)}(\lambda)}{(m-1)!} \\
     0 & \ddots &  \ddots  & \ddots  &\vdots  \\
     \vdots    &\ddots  & \ddots & \ddots &    
       \frac{f''(\lambda)}{2}   \\
    \vdots & &\ddots  & f(\lambda) & f'(\lambda)     \\ 
   0 & \cdots &\cdots  & 0 & f(\lambda)   
\end{pmatrix}
$$
 Above,
 $f^{(m-1)}(\lambda)$ denotes the $(m-1)$st derivative of $f$ at the point $\lambda$. In particular, the first row consists of  first coefficients of the Taylor extension of   $f$ 
at $\lambda$.  
Combining this formula with the definition of $D$ we see that (in   dimension $m$ and thinking that $\lambda$ is a function)
$$
\left(\frac{\partial}{\partial x^1}, \dots , \frac{\partial}{\partial x^m}\right) f(J_m(\lambda))= 
\left(\frac{\partial \lambda}{\partial x^1}, \dots , \frac{\partial \lambda}{\partial x^m}\right) f'(J_m(\lambda)), 
$$
where $f'(t)$ is the $t$-derivative of the polynomial $f(t)$.

Therefore, 
$
D(K^n) - \sigma_1 D(K^{n-1})-...-\sigma_{n-1} D(K)$ is equal to $$
\left(\frac{\partial \lambda}{\partial x^1}, \dots , \frac{\partial \lambda}{\partial x^m}\right) \left( \frac{d}{dt}(\chi_K(t))_{|t=K}\right).$$
Using that $$\frac{d}{dt} (\chi_K(t)) = \sum_s m_s (t- \lambda_s)^{m_s-1}  \prod_{r\ne s} (t- \lambda_r)^{m_r}
$$
we obtain  that $
D(K^n) - \sigma_1 D(K^{n-1})-...-\sigma_{n-1} D(K)$ is equal to 
\begin{eqnarray*}
    \left(\frac{\partial \lambda}{\partial x^1}, \dots , \frac{\partial \lambda}{\partial x^m}\right) &&  
\Big( m_1 \left(A-\lambda_1 \operatorname{Id}\right)^{m_1-1} \left(A-\lambda_2 \operatorname{Id}\right)^{m_2}\cdots  \left(A-\lambda_2 \operatorname{Id}\right)^{m_k}\\ 
& +&  m_2 \left(A-\lambda_1 \operatorname{Id}\right)^{m_1}
\left(A-\lambda_2 \operatorname{Id}\right)^{m_2-1}\left(A-\lambda_3 \operatorname{Id}\right)^{m_3}\cdots  \left(A-\lambda_2 \operatorname{Id}\right)^{m_k} \\  & + & \cdots  \\ &+& 
m_k \left(A-\lambda_1 \operatorname{Id}\right)^{m_1}
\cdots \left(A-\lambda_{k-1} \operatorname{Id}\right)^{m_{k-1}} \left(A-\lambda_2 \operatorname{Id}\right)^{m_k-1}\Big).
\end{eqnarray*}
Next, we see that for the first $m_1$ components 
of $
D(K^n) - \sigma_1 D(K^{n-1})-...-\sigma_{n-1} D(K)$
only the first product in the sum above is relevant, and it gives zeros on the first $m_1- 1$ places and $$m_1 \frac{\partial \lambda_1}{ \partial x^1 } 
 \prod_{s\ne 1} (\lambda_1- \lambda_s)^{m_s}$$ on the place number $m_1$. Similarly, for the components from $m_1+1$ to $m_1+m_2 $ only the second product in the sum above is relevant. Then, the components  from 
 $m_1+1$ to $m_1+m_2-1 $ are zero and the component number $m_1+m_2$ equals  $$ m_2 \frac{\partial \lambda_2}{ \partial x^{m_1+1} } 
 \prod_{s\ne 2} (\lambda_2- \lambda_s)^{m_s} $$ and so on.  Lemma \ref{lem:1}, and therefore Theorem \ref{thm:3}, are proved. 

\section{Conclusion and outlook}
Results of our paper give the following test   to understand whether a given $(1,1)$-tensor field $K^i_j$  generates  separation of variables  for a certain (a priori unknown) metric: one needs  to  answer the following four test questions. 

\begin{itemize}
\item[Q1]  Is  the discriminant of the characteristic polynomial of $K$ is different from zero? 

\item[Q2] Is the Haantjes torsion of $K$    zero? 

\item[Q3] Is  \eqref{eq:M3} fulfilled?

\item[Q4] Is the  $(1,3)$-tensor  $P$ given by \eqref{eq:last}   identically zero? 
\end{itemize}

If all questions are  answered by ``Yes" then there (locally) exists a metric  for which $K$ generates separation of variables. If at least one of the answers is ``No", then  no such metric exists. 

The  related calculations can  be made in any coordinate system. The discriminant, Haantjes torsion,  \eqref{eq:M3} and \eqref{eq:last} are given by explicit formulas  which are polynomial expressions in the components of $K$ and their first and second derivatives. They    can be implemented in standard  computer algebra packages. 

Of course, provided the tensor field  $K$ passed all the tests above, it could be 
interesting to find a metric $g$ for which $K$ generates separation of variables.  From the proof of Theorem \ref{thm:2} we know that in ``diagonal'' coordinates for $K$ one can find such a  metric by integrating certain closed 1-forms explicitly constructed by $K$ (and the metric is defined up to an arbitrary choice of $n$ functions of one variable). 
We do not have such a nice answer in an arbitrary coordinate system. Let us note though that the equations for  the components  of the   metric $g$  are linear first order equations. Indeed, if the components $K^i_j$ are  given and $g_{ij}$ are viewed as unknown functions\footnote{We essentially have $\frac{n(n+1)}{2}$ unknown functions since $g_{ij}$ must satisfy the relation $K^s_jg_{si}=K^s_ig_{sj} $}, 
the first term of \eqref{eq:killingequation} reads 
$$
\begin{array}{rl}K_{ij,k}= g_{si} K^s_{j,k}= &  g_{si} \frac{\partial K^s_j}{\partial x^k} - g_{si}  K^s_r \Gamma^r_{jk}+ g_{si} K^r_j \Gamma^s_{rk}\\
=& g_{si} \frac{\partial K^s_j}{\partial x^k} -   K^s_i \Gamma_{ jk,s  }+  K^r_j \Gamma_{rk,i}.
\end{array}
$$
Since the Christoffel symbols of  the  first kind, $$
\Gamma_{ij,k}:= \frac{1}{2} (g_{ik,j}+g_{jk,i}-g_{ij,k}),$$
are clearly linear in $g_{ij}$, we see that  $K_{ij,k}$ is given by a first order  linear expression in $g$ with coefficients constructed by $K^i_j$. Similarly, the other two terms in \eqref{eq:killingequation}  are linear in $g_{ij}$ so the system 
\eqref{eq:killingequation} is a linear (over-determined) system of the first order on the components of $g$. Theorems \ref{thm:3} and \ref{thm:4}  can be viewed as a  geometric  way to write  the compatibility conditions for this system, under the additional assumption that 
the Haantjes torsion of $K$ vanishes and $K$ has $n$ distinct eigenvalues. It is  remarkable  though that, differently from many visually  similar problems, see 
e.g. \cite{K1,K2},  compatibility conditions can be written explicitly, in a closed and relatively simple form; moreover, no branching appears.

Our motivation to consider  the $(1,1)$ tensor $K$  as a ``main''  object is related to the research programme ``Nijenhuis Geometry''  suggested in \cite{BMMT,Nijenhuis1}. Our further goals within this programme which may use the results of the present paper will include understanding of the natural analogue of 
semi-hamiltonicity condition  and of separation of variables  in the case when the corresponding $(1,1)$ tensor is gl-regular but still  has Jordan blocks. First results in this direction were already obtained in \cite{BKM3,BKM2}.

\subsubsection*{ Acknowledgements.} Vladimir Matveev  thanks the DFG (projects 455806247 and  529233771)  and  ARC Discovery Programme   (grant DP210100951) for the support. Most   results were obtained during  research visits  of VM to the University of New South Wales and Sydney University  supported by the  Sydney Mathematics Research Institute; VM thanks the UNSW,  SMRI and Sydney University  for their hospitality. We thank A. Bolsinov and  E. Ferapontov for useful discussions.

\printbibliography

\end{document}